\documentclass[12pt]{amsart}
\usepackage{epsfig}
\usepackage{amsmath}
\usepackage{amsfonts}
\usepackage{amssymb}
\usepackage{epigraph, fancyhdr}

\def\A{\mathcal A}

\def\D{\mathcal D}
\def\F{\mathcal F}
\def\K{\mathcal K}

\def\J{\mathcal J}
\def\L{\mathcal L} 

\def\O{\mathcal O}
\def\T{\mathcal T}

\def\1{\mathbf 1}
\def\M{{\overline{\mathcal M}}}

\def\QQ{\mathbb Q}
\def\ZZ{\mathbb Z}

\def\Res{\operatorname{Res}}
\def\Ind{\operatorname{Ind}}

\def\hat{\widehat}
\def\tilde{\widetilde}
\def\p{\partial}
\def\a{\alpha}
\def\b{\beta}

\def\f{{\mathbf f}}
\def\g{{\mathbf g}}
\def\t{{\mathbf t}}
\def\x{{\mathbf x}}
\def\y{{\mathbf y}}
\def\v{{\mathbf v}}

\def\gs{\sigma}
\def\gL{\Lambda}

\def\h{\hbar}

\def\lan{\langle}
\def\ran{\rangle}

\def\ft{\operatorname{ft}}
\def\ct{\operatorname{ct}}
\def\td{\operatorname{td}}
\def\ch{\operatorname{ch}}

\def\fake{\operatorname{fake}}
\def\tr{\operatorname{tr}}

\def\z{\zeta}

\renewcommand{\Delta}{\triangle}
\makeatletter
\DeclareFontFamily{OMX}{MnSymbolE}{}
\DeclareSymbolFont{MnLargeSymbols}{OMX}{MnSymbolE}{m}{n}
\SetSymbolFont{MnLargeSymbols}{bold}{OMX}{MnSymbolE}{b}{n}
\DeclareFontShape{OMX}{MnSymbolE}{m}{n}{
    <-6>  MnSymbolE5
   <6-7>  MnSymbolE6
   <7-8>  MnSymbolE7
   <8-9>  MnSymbolE8
   <9-10> MnSymbolE9
  <10-12> MnSymbolE10
  <12->   MnSymbolE12
}{}
\DeclareFontShape{OMX}{MnSymbolE}{b}{n}{
    <-6>  MnSymbolE-Bold5
   <6-7>  MnSymbolE-Bold6
   <7-8>  MnSymbolE-Bold7
   <8-9>  MnSymbolE-Bold8
   <9-10> MnSymbolE-Bold9
  <10-12> MnSymbolE-Bold10
  <12->   MnSymbolE-Bold12
}{}

\let\llangle\@undefined
\let\rrangle\@undefined
\DeclareMathDelimiter{\llan}{\mathopen}%
                     {MnLargeSymbols}{'164}{MnLargeSymbols}{'164}
\DeclareMathDelimiter{\rran}{\mathclose}%
                     {MnLargeSymbols}{'171}{MnLargeSymbols}{'171}
\makeatother
\title[General theory]
      {Permutation-equivariant \\ quantum K-theory VII. \\
      General theory}

\author[A. Givental]{Alexander GIVENTAL}
\thanks{This material is based upon work supported by the National 
Science Foundation under Grant DMS-1007164, and by the IBS Center for Geometry 
and Physics, POSTECH, Korea.} 

\date{August 4, 2015}

\begin{document}

\begin{abstract}

  We introduce K-theoretic GW-invariants of mixed nature: permutation-equivariant in some of the inputs and ordinary in the others, and prove the ancestor-descendant correspondence formula. In genus 0, combining this with adelic characterization, we derive that the range $\L_X$ of the big J-function in permutation-equivariant theory is overruled.
  %and that the ruling spaces are $D_q$-modules in Novikov's variables. We describe a large group of symmetries of $L_X$ defined in terms of $q$-differenc operators. We use it to explicit reconstruction of $L_X$ from one point. We apply the result to toric $X$, when such a point is given by the $q$-hypergeometric function.         

\end{abstract}

\maketitle

\section*{The string and dilaton equations}

We return to the introductory setup of Part I, and introduce {\em mixed genus-$g$ descendant potentials} of a compact K\"ahler manifold $X$:
\[ \F_g(\x, \t):=\sum_{k\geq 0,n\geq 0, d} \frac{Q^d}{k!}\lan \x(L),\dots,\x(L); \t(L),\dots, \t(L)\ran_{g,k+n,d}^{S_n}.\]
The first $k$ seats are occupied by the input $\x =\sum_{r\in \ZZ} x_r q^r$, which is a Laurent polynomial in $q$ with vector coefficients $x_r \in K^0(X)\otimes \gL$. We assume that $\gL$ includes Novikov's variables as well.
The last $n$ seats are occupied by similar inputs
$\t=\sum_{r\in \ZZ} \t_rq^r$, $\t_r\in K^0(X)\otimes \gL$, and only these inputs are considered {\em permutable} by renumberings of the marked points. Most of the time we will assume that $\t(q)=t$ is constant in $q$, i.e. that the permutable inputs do not involve the cotangent line bundles $L_i$. 

We will first treat these generating functions as objects of the {\em ordinary}, i.e.  permutation-{\em non}-equivariant, quantum K-theory, depending however on the parameter $t$.
Our nearest aim is to extend to this family of theories some basic facts from the ordinary GW-theory, starting with the genus-0 string and dilaton equations.

On the moduli space $X_{g,m+1,d}$, along with the line bundles $L_i$ formed by the cotangent lines to the curves at the $i$th marked point, consider the line bundles $\tilde{L}_i:=\ft_1^*(L_i)$, $i\geq 2$, where $\ft_1: X_{g, 1+m,d} \to X_{g,m,d}$ is the map defined by forgetting the first marked point. In the {\em genus 0} case, it is clear that
\[ \lan 1, \x_1(\tilde{L}),\dots,\x_k(\tilde{L}); t,\dots \ran_{0,1+k+n,d}^{S_n} = 
\lan \x_1(L), \dots, \x_k(L) ; t,\dots \ran_{0,k+n,d}^{S_n}.\]
On the other hand, it is well-known how to compare $L_i$ and $\tilde{L}_i$.
The fibers of these line bundles coincide everywhere outside the section
$\gs_i: X_{g,m,d}\to X_{g,1+m,d}$ defined by the $i$th marked point, while $\gs_i^*(\tilde{L}_i)=L_i$, and $\gs_i^*(L_i)=1$. In other words,
$1-\tilde{L}_i/L_i = (\gs_i)_*1$, and hence $\tilde{L}_i=L_i-(\gs_i)_*1$.
Taking into account that $L_i ((\gs_i)_*1)=(\gs_i)_*1$, and that $((-\gs_i)_*1)^r=(-\gs_i)_*(L_i-1)^{r-1}$,  we find by Taylor's formula (and omitting the subscript $i$):
\begin{align*} \x(\tilde{L})- \x(L)&=\sum_{r>0} \frac{\x^{(r)}(L)}{r!}(-\gs_*1)^r= \\ 
  -\gs_*\left( \sum_{r>0} \right. & \left.\frac{\x^{(r)}(1)}{r!}(L-1)^{r-1} \right) =  -\gs_* \frac{\x(L)-\x(1)}{L-1} .\end{align*}
Note that the divisors $\gs_i$ for different $i$ are disjoint, and that $\gs_i^*(L_j)=L_j$ if $j\neq i$.  Thus
\begin{align*} \lan 1, \x_1(L),\dots, \x_k(L); t,\dots\ran_{0,1+k+n,d}^{S_n} = \lan &\x_1(L),\dots, \x_k(L); t,\dots \ran_{0,k+n,d} +\\ 
\sum_{i=1}^k\lan \dots, \x_{i-1}(L), \frac{\x_i(L)-\x_i(1)}{L-1}, &\x_{i+1}(L),\dots ; t, \dots, t \ran_{0,k+n,d}^{S_n}.\end{align*}
This computation is quite standard, since it does not interfere with the permutable inputs, as long as those don't contain line bundles $L_i$. 

\medskip

{\tt Proposition 1} (string equation). {\em Let $V$ be the linear vector field on the space of vector-valued Laurent polynomials in $q$ defined by
  \[ V (\y) := \frac{\y(q)-\y(1)}{1-q}.\]
In the genus-0 descendent potential $F_0(\x,t)$, introduce the {\em dilaton shift} of the origin: $\y(q)=1-q+t+\x(q)$. Then
  \[ L_V (\F_0(\y+q-1-t), t)) = \F_0 (\y+q-1-t), t)+\frac{(\y(1),\y(1))}{2}-\left(\frac{\Psi^2(t)}{2},1\right),\]
  where $(a,b):=\chi (X; ab)$ is the $\gL$-valued K-theoretic Poincar{\'e} pairing, and $\Psi^2$ is the 2nd Adams operation on $K^0(X)\otimes \gL$.} 
  
\medskip

{\tt Proof.} The linear vector field $V$ becomes inhomogeneous in the unshifted coordinate system:
\[ \frac{\y(q)-\y(1)}{1-q}=\frac{\x(q)-\x(1)}{1-q}+1.\]
Applying the previous, down-to-earth form of the string equation to
\[ \F_0 (\x) := \sum_{k,n,d}\frac{Q^d}{k!}\lan \x(L),\dots,\x(L);t,\dots,t\ran_{0,1+k+n,d}^{S_n},\]
we gather that
\[ L_V (\F_0(\x,t)) = \F_0(\x) + \text{terms $\lan 1, \dots\ran^{S_n}_{0,3,0}$ with $d=0$ and $k+n=2$}.\]
Since $X_{0,3,0}=X\times \M_{0,3}=X$, and $L=1$ on $\M_{0,3}$, these terms are
\[ \frac{1}{2}(\x(1),\x(1))+(\x(1),t)+\frac{1}{2}(t,t)/2-
\frac{1}{2}(\Psi^2(t),1).\]
The last two terms come from
\[ \lan 1; t,t\ran_{0,3,0}^{S_2}=\frac{1}{|S_2|}\sum_{h\in S_2} \tr_h (t^{\otimes 2}).\]
All but the last one add up to $(y(1),y(1))/2$. \qed

\medskip

Consider now correlators
\[ \lan L-1, \x(L),\dots,\x(L); \t(L), \dots, \t(L)\ran_{0,1+k+n,d}^{S_n}.\]
The line bundle $L_1$ over $X_{g,1+k+n,d}$ differs from the dualizing sheaf to the fibers of the forgetting map $\ft_1:X_{g,1+k+n,d}\to X_{g,k+n,d}$ by the divisor of the marked points. The spaces $H^0(\Sigma, L)$ are formed by 
holomorphic differentials on $\Sigma$ with at most 1st order poles at the markings, and with at most 1st order poles at the nodes with zero residue sum at each node. In genus $0$, if $k+n>0$, then $H^1(\Sigma, L-1)=0$, while the holomorphic differentials are uniquely determined by the residues at the marked points subject to the constrains that the total sum is $0$. The residues {\em per se} form trivial bundles, but those at the permutable marked points form the standard Coxeter representation of $S_n$, induced from the trivial representation of $S_{n-1}$. Thus,
\[  (ft_1)_*(L-1)=k-2+\Ind_{S_{n-1}}^{S_n}(1),\]
and this answer is correct even when $k=n=0$ (in which case $H^1(\Sigma, L)=H^0(\Sigma, 1)^*=1$).
On the other hand, $L-1$ vanishes on the sections $\gs_i:X_{g,k+l,d}\to X_{g,1+k+n,d}$ defined by the markings, where the differences between $\tilde{L_i} - L_i$, $i>1$, are supported. We find that
\begin{align*} \lan L-1, \dots,\x(L); \t(L), \dots \ran_{0,1+k+n,d}^{S_n} = (k-2) \lan \dots,\x(L);\t(L),\dots \ran_{0,k+n,d}^{S_n} \\
+\lan \x(L) \dots,\x(L), \t(L); \t(L),\dots,\t(L)\ran_{0,k+n,d}^{S_{n-1}}.\end{align*}
We use here that for any $S_n$-module $V$,
\[ \left( V\otimes \Ind_{S_{n-1}}^{S_n}(1)\right)^{S_n}=\left( \Res_{S_{n-1}}^{S_n}(V)\right)^{S_{n-1}}.\]

{\tt Proposition 2} (dilaton equation). {\em The genus-0 descendent potential $\F_0$ in dilaton-shifted coordinates satisfies the following homogeneity condition:
  \[ L_E (\F_0(\y+q-1-\t,\t) = 2 \F_0(\y+q-1-\t,\t) - (\Psi^2(\t(1)),1),\]
where $E$ is the Euler vector field $E(\y)=\y$ in the linear space of vector-valued Laurent polynomials $\y(q)$.}

\medskip

{\tt Proof.} The exceptional terms
\[ \frac{1}{2}\lan L-1, \x(L),\x(L)\ran_{0,3,0}+\lan L-1, \x(L); \t(L)\ran_{0,3,0}^{S_1}+\lan L-1;\t(L),\t(L)\ran_{0,3,0}^{S_2}\]
all vanish except for the trace of the non-trivial element in $S_2$, which acts by $-1$ one the cotangent line $L$. This makes $L-1$ on $\M_{0,3}$ equal to $-2$ (rather than $0$), and results in the constant $-(\Psi^2(\t(1)),1)$. Therefore the identity derived above yields:
\begin{align*} \sum_{k,n,d}\frac{Q^d}{k!}\lan 1-L+\t(L)+\x(L), \x(L),\dots,\x(L);
  \t(L),\dots,\t(L)\ran_{0,1+k+n,d}^{S_n} \\
  = 2\sum_{k,n,d}\frac{Q^d}{k!}\lan \x(L),\dots,\x(L);\t(L),\dots,\t(L)\ran_{0,k+n,d}^{S_n} - (\Psi^2(\t(1)),1),\end{align*}
  which after shift $\y(q):=1-q+\t(q)+\x(q)$ becomes what we claimed.  \qed
  
{\tt Remark.} Note that we have proved this allowing the permutable input $\t$,  i.e. the {\em parameter} of $\F_0$ to depend on $q$. 

\section*{A WDVV-equation}

Let us introduce the gadget
\[ \llan A_1, \dots, A_m \rran_{g,m}:=\sum_{l,n,d}\frac{Q^d}{l!}\lan A_1,\dots,A_m;
\tau, \dots, \tau; t, \dots, t\ran_{g,m+l+n,d}^{S_n},\]
for the generating function of $\tau, t \in K^0(X)\otimes \gL$, and the meaning of the inputs $A_i$ to be specified.

\medskip

Along with the {\em Poincar\'e metric} $g_{\a\b}=(\phi_{\a},\phi_{\b})$ on $K^0(X)$, where $\{ \phi_{\a} \}$ is a basis, introduce the non-constant metric
\[ G_{\a\b}:= g_{\a\b}+\llan \phi_{\a},\phi_{\b}\rran_{0,2}.\]
Note that the inverse tensor has the form
\begin{align*} G^{\a\b}=g^{\a\b}-&\llan\phi^{\a},\phi^{\b}\rran_{0,2}+\sum_{\mu}\llan \phi^{\a},\phi^{\mu}\rran_{0,2}\llan\phi_{\mu},\phi^{\b}\rran_{0,2}\\ & -\sum_{\mu,\nu}
\llan\phi^{\a},\phi^{\mu}\rran_{0,2}\llan\phi_{\mu},\phi^{\nu}\rran_{0,2}\llan\phi_{\nu},\phi^{\b}\rran_{0,2}+\dots,\end{align*}
where $\{\phi^{\a}\}$ is the basis Poincar{'e}-dual to $\{\phi_{\a}\}$. 

\medskip
  
{\tt Proposition 3} (WDVV-equation). {\em For all $\phi,\psi \in K^0(X)\otimes \gL$, 
\begin{align*} (\phi, \psi)&+(1-xy)\llan  \frac{\phi}{1-xL},\frac{\psi}{1-yL}\rran_{0,2} = \\   
\sum_{\a,\b} &\left( (\phi,\phi_{\a})+\llan \frac{\phi}{1-xL},  \phi_{\a}\rran_{0,2}\right) \,G^{\a\b}\left( (\phi_{\b},\psi)+\llan  \phi_{\b},\frac{\psi}{1-yL}\rran_{0,2}\right) .\end{align*} }

\medskip

{\tt Proof.} The standard WDVV-argument consists in mapping moduli spaces of genus-0 stable maps with $4+$ marked points to the Deligne-Mumford space $\M_{0,4}$, and considering the inverse image of a typical point, i.e., in other words, fixing the cross-ratio of the first $4$ marked points. When the cross-ratio degenerates into one of the special values $0,1,\infty$, the curves become reducible, with the 4 marked points split into pairs between the two glued pieces in 3 different ways. The WDVV-equation expresses the equality between the three gluings. 

We apply the argument to the inputs of the 4 marked points equal to $1, 1, \phi/(1-xL)$, and $\phi/(1-yL)$, and arrive at the following identity (see Figure 1):
\begin{align*} \sum_{\a,\b} \llan 1, \frac{\phi}{1-xL}, \phi_{\a}\rran_{0,3}G^{\a\b}\llan \phi_{\b}, \frac{\psi}{1-yL},1 \rran_{0,3} = \\
    \sum_{\a,\b}\llan 1 , 1, \phi_{\a} \rran_{0,3}G^{\a\b}\llan \phi_{\b}, \frac{\phi}{1-xL}, \frac{\psi}{1-yL} \rran_{0,3}.\end{align*} 

\begin{figure}[htb]
\begin{center}
\epsfig{file=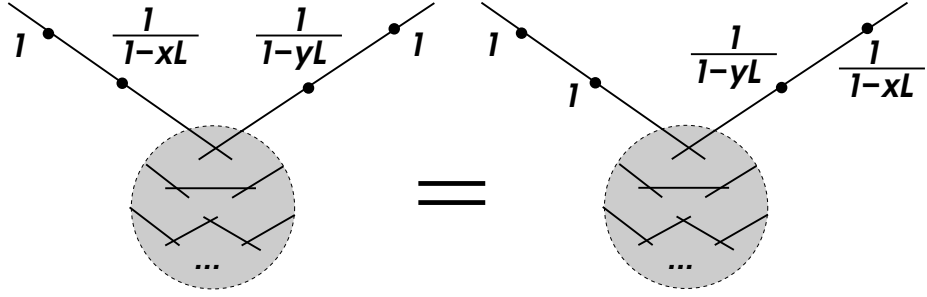} 
\end{center}
\caption{WDVV equation}
\end{figure}

As it is explained in \cite{GiK}, in K-theory the WDVV-argument encounters the following subtlety. The virtual divisor obtained by fixing the cross-ratio and passing to any of the three limits, has self-intersections, represented by curves with more than 2 components (as shown on Figure 1 in shaded areas). As a result, the structure sheaf of the divisor before the limit is identified with the alternated sum of the structure sheaves of all the self-intersection strata on a manner of the exclusion-inclusion formula. In the identity, this is taken care of by the pairing which involves the tensor $G^{\a\b}$.

It only remains to apply the string equation. Since
\[ \frac{1}{L-1}\left(\frac{1}{1-qL}-\frac{1}{1-q}\right)=\frac{x}{(1-x)}\frac{1}{(1-qL)},\]
and since $L=1$ on $X_{0,3,0}=X\times \M_{0,3}=X$, we have 
\begin{align*} \llan 1, \frac{\phi}{1-qL}, \phi_{\a}\rran_{0,3} &= \frac{(\phi, \phi_{\a})}{1-q} +\left(1+\frac{q}{1-q}\right)\llan \frac{\phi}{1-qL},  \phi_{\a}\rran_{0,2} ;\\
  \sum_{\a,\b}\llan 1 , 1, \phi_{\a} \rran_{0,3}G^{\a\b}&\llan \phi_{\b}, \frac{\phi}{1-xL}, \frac{\psi}{1-yL} \rran_{0,3} = 
\llan 1, \frac{\phi}{1-xL}, \frac{\psi}{1-yL} \rran_{0,3} \\ = \frac{(\phi, \psi)}{(1-x)(1-y)}&+
\left(1+\frac{x}{1-x}+\frac{y}{1-y}\right) \llan \frac{\phi}{1-xL}, \frac{\psi}{1-yL}\rran_{0,2}.\end{align*}
The result follows. \qed 

\section*{The loop space formalism}

Here we interpret the string, dilaton, and WDVV-equations using symplectic linear algebra in the space $\K$ of rational functions of $q$ with vector  coefficients from $K^0(X)\otimes \gL$. To be more precise, we assume that
elements of $\K$ are such rational functions {\em modulo any power of Novikov's variables} (or in any other topology that may turn out useful in future). We equip $\K$ with symplectic form
\[ \Omega (\f,\g) := -\Res_{q=0,\infty} (\f(q^{-1}), \g(q))\, \frac{dq}{q}.\]
We identify $\K$ with $T^*\K_{+}$ where $\K_{+} \subset \K$ is the Lagrangian subspace consisting of vector-valued Laurent polynomials in $q$ (in the aforementioned topological sense) by picking the complementary Lagrangian subspace $\K_{-}$ consisting of rational functions of $q$ regular at $q=1$ and
vanishing at $q=\infty$. We encode K-theoretic genus-0
GW-invariant of $X$ by the {\em big J-function}
\begin{align*} \J(\x,\t)&:=1-q+\t(q)+\x(q)+\\
  &\sum_{\a,k,n,d}\phi^{\a}\frac{Q^d}{k!}\lan \frac{\phi_{\a}}{1-qL}, \x(L),\dots,\x(L);\t(L),\dots,\t(L)\ran_{0,1+k+n,d}^{S_n}.\end{align*}

{\tt Proposition 4.} {\em The big J-function is the dilaton-shifted graph of the differential of the genus-0 descendent potential $\F_0$:
\[ \J(\x,\t)= 1-q+\t+\x+d_{\x}\F_0(\x,\t).\]}

\medskip

{\tt Proof.} For every $\v \in \K_{+}$, we have: 
\begin{align*} L_{\v}\F_0&=\sum_{k,n,d}\frac{Q^d}{k!}\lan \v(L),\x(L),\dots,\x(L);\t(L),\dots,\t(L)\ran_{0,1+k+n,d}^{S_n} =\\
  \sum_{k,n,d}\frac{Q^d}{k!}& \lan \Res_{q=L} \frac{\v(q)}{(1-L/q)}\, \frac{dq}{q}, \x(L),\dots,\x(L);\t(L),\dots,\t(L)\ran_{0,1+k+n,d}^{S_n} \\
  &=-\Res_{q=0,\infty} (\J(q^{-1}), \v(q))\, \frac{dq}{q} = \Omega (\J, \v).\end{align*}

{\tt Corollary 1} (dilaton equation). {\em For a fixed value of the parameter $\t$, the range of the J-function $\x \mapsto \J (\x,\t)$ is a Lagrangian cone $\L_{\t} \subset \K$ with the vertex at the origin.}

\medskip

{\tt Proof.} Differentiating the dilaton equation for $\F_0$, we find that 1st derivatives of $\F_0$ are homogeneous of degree 1. \qed  

\medskip

{\tt Corollary 2} (string equation). {\em For any $t\in K^0(X)\otimes \gL$, the linear vector field $\f \mapsto \f/(1-q)$ on $\K$ is tangent to $\L_t$.}

\medskip

{\tt Remark.} We will see later that this is true for any $\L_{\t}$, and not only for $\t$ independent of $q$. 

\medskip

{\tt Proof.}  Subtracting from the string equation for $\F_0$ derived in the previous section a half of the dilaton equation for $\F_0$, we obtain a Hamilton-Jacobi equation $L_{V-E/2} \F_0 = (\y(1),\y(1))/2$. It expresses the fact that the quadratic Hamiltonian corresponding to this equation vanishes on $\L_t$, and hence the Hamiltonian vector field is tangent to $\L_t$. We will show that this Hamiltonian vector field is
\[ W\f :=\frac{\f}{1-q}-\frac{\f}{2}.\]
Due to Corollary 1, $\f\mapsto \f/2$ is tangent to $\L_t$, and the result about
$\f\mapsto \f/(1-q)$ would follow. 

The hamiltonian of $W$ is $H(\f):=\Omega (\f,W\f)/2 =\Omega (\f, \f/(1-q))/2$. Using the projections $\f_{\pm}$ of $\f\in \K$ to $\K_{\pm}$, we compute $2H(\f)$:
\begin{align*} \Omega \left(\f, \frac{\f}{1-q}\right) &= \Omega \left(\f_{+}+\f_{-}, \frac{\f_{+}(1)}{1-q}+
  \frac{\f_{+}-\f_{+}(1)}{1-q}+\frac{\f_{-}}{1-q} \right)=\\
-\Omega\left(\frac{\f_{+}(1)}{1-q},\f_{+}\right)&+\Omega\left(\f_{-},\frac{\f_{+}-\f_{+}(1)}{1-q}\right)+\Omega \left(\frac{\f_{+}}{1-q^{-1}},\f_{-}\right) = \\
  \Res_{q=0,\infty}&\left(\frac{\f_{+}(1)}{1-q^{-1}},\f_{+}(q)\right)\frac{dq}{q} +\Omega\left(\f_{-},\frac{\f_{+}-2\f_{+}(1)+q\f_{+}}{1-q}\right)+\\
  \Omega \left(\f_{-},\frac{\f_{+}(1)}{1-q}\right)&=
 -(\f_{+}(1),\f_{+}(1))+2\, \Omega \left(\f_{-},\frac{\f_{+}-\f_{+}(1)}{1-q}-\frac{\f_{+}}{2}\right)+0.\end{align*}
The last non-zero term is twice the Hamilton function of the vector field
$\y \mapsto \frac{\y(q)-\y(1)}{1-q}-\y(q)/2$ on $\K_{+}$, i.e. $V-E/2$, lifted to the cotangent bundle in the standard way. The first non-zero term is twice
$-(\y(1),\y(1))/2$. Thus, the quadratic hamiltonian $H$ is exactly as claimed.
\qed

\medskip

Introduce the operator $S: K^0(X)\otimes \gL \to \K_{-}$ defined by
\[ S (q)\, \phi = \sum_{\a,\b}\left((\phi,\phi_{\a})+\llan \frac{\phi}{1-L/q}, \phi_{\a} \rran_{0,2} \right) G^{\a\b}\phi_{\b},\]
The operator depends on the parameter $\tau \in K^0(X)\otimes \gL$. For each value of the parameter, it can be considered as an operator-valued rational function of $q$ (a ``loop group'' element), and in this capacity extends to
 a map $S: \K\to \K$ commuting with multiplications by scalar rational functions of $q$. The WDVV-identity from the previous section can be written as
 \[ (1-xy)\llan \frac{\phi}{1-xL}, \frac{\psi}{1-yL}\rran_{0,2}  = (\phi,\psi)+\left(S^*(y^{-1})S(x^{-1})\phi, \psi \right), \]
 where
 \[ S^*(q)\psi = \psi +\sum_{\a\b}\llan \psi, \frac{\phi_{\a}}{1-L/q}\rran_{0,2}g^{\a\b}\phi_{\b}\]
   is the operator adjoint to $S(q)$ with respect to the inner product $(g_{\a\b})$ on the domain space, and $(G_{\a\b})$ on the target space.
It follows that
\[ S^*(q^{-1})S(q)=1,\ \ \text{and hence}\ \ S(q)S^*(q^{-1})=1.\]
This means that $S: (\K ,\Omega) \to (\K, \bar{\Omega})$ provides a symplectic isomorphism between two symplectic structures on the loop space: $\Omega$, based on the metric tensor $(g_{\a\b})$, and $\bar{\Omega}$, based on the metric tensor $(G_{\a\b})$ (and depending therefore on the parameter $\tau \in K^0(X)\otimes \gL$). The inverse isomorphism is given by
\[ S^{-1}(q)=S^*(q^{-1}).\]    
Furthermore, the {\em quantizations} $\hat{S}$ and $\widehat{S^{-1}}$ provide isomophisms between the corresponding Fock spaces, which in their formal version consist of expressions
\[  \D (\y) = e^{\textstyle \F_0(\y)/\h + \F_1(\y)+\h \F_2(\y)+ \h^2 \F_3(\y)+\cdots },\ \ \ \y\in \K_{+},\]
where $\F_g$ is a sequence of scalar-valued functions on $\K_{+}$.

\medskip

{\tt Proposition 5.} {\em The action of the quantized operator $S^{-1}$ on an element $\A$ of the Fock space, corresponding to the symplectic form $\bar{\Omega}$, is given by
  \[ (\widehat{S^{-1}} \A) (\y) = e^{\llan \y(L), \y(L)\rran_{0,2}/\h}\, \A ([S(q) \y(q)]_{+}),\]
    where $[f(q)]_{+}$ denotes taking the Laurent polynomial part of rational function $\f$, i.e. the projection $\K\to \K_{+}$ along $\K_{-}$.}

\medskip

{\tt Proof.}   Generally speaking, quantization of linear symplectic transformations $T$ is defined as $\exp \widehat{\ln T}$, where $\widehat{\ln T}$ is quantization of the quadratic hamiltonian according to the standard rules \cite{GiQ}. Namely, in Darboux coordinates on $\K=T^*\K_{+}$, 
\[ \widehat{q_{\a}q_{\b}} = \h^{-1}q_{\a}q_{\b}, \  \widehat{q_{\a}p_{\b}}=q_{\a}\p_{q_{\b}}, \ \widehat{p_{\a}p_{\b}}=\h \p_{p_{\a}}\p_{p_{\b}}.\]
The operators $S^{\pm 1}$ have the form of the composition of the operator $G^{\pm 1}$ identifying the metric: $(G\phi_{\mu},\phi_{\nu})=G_{\mu\nu}$,
and the operator $(\K,\Omega) \to (\K,\Omega)$ which is the identity modulo
$\K_{-}$. This means that the quadratic hamiltonian of $\ln S^{-1}G$ contains only $pq$-terms and $q^2$-terms, but no $p^2$-terms. Therefore $\widehat{S^{-1}}:=\exp(\widehat{\ln S^{-1}G}) G^{-1}$ will act by a linear change of variables followed by the multiplication by a quadratic form, both depending on $S$. The answer in the finite form is given by Proposition 5.3
in \cite{GiQ}:
\[ (\widehat{S^{-1}} \A) (\y) = e^{W(\y,\y)/\h}\, \A ([S(q) \y(q)]_{+}),\]
where the symmetric bilinear form $W$ is determined by
\[ W(\x,\y) = (\Omega \otimes \Omega) \left(\frac{S^*(x^{-1})S(y^{-1})-1}{1-xy}, \x(x)\otimes \y(y)\right).\]
Using the WDVV-equation, we find
\begin{align*} W(\x,\y) =&
  \Res_{x={0,\infty}} \Res_{y=0,\infty} \left\llan \frac{\x(x)}{1-L/x}, \frac{\y(y)}{1-L/y}\right\rran_{0,2} \frac{dx}{x}\frac{dy}{y} \\
 = &\llan \x(L),\y(L)\rran_{0,2}.\end{align*}

\section*{Ancestor -- descendent correspondence}  
 
We introduce {\em ancestor potentials}
\begin{align*} \bar{\F}_g(\x, \tau, t)&:=\sum_{k\geq 0,d} \frac{1}{k!}\llan \x(\bar{L},\dots,\x(\bar{L})\rran_{g,k}= \\ \sum_{k,l,n,d}& \frac{Q^d}{k!l!}\lan \x(\bar{L}),\dots,\x(\bar{L}); \tau,\dots,\tau; t,\dots, t\ran_{g,k+l+n,d}^{S_n}, \end{align*}
where $\tau, t \in K^0(X)\otimes \gL$, and $\bar{L}$ in the $i$th position of the correlator represents the line bundle $\bar{L}_i$ over the moduli space $X_{g,k+l+n,d}$ of stable maps to $X$, obtained by pulling back the universal cotangent line bundle at the $i$th marked point over the Deligne-Mumford space $\M_{g,k}$ by the {\em contraction} map $\ct: X_{g,k+l+n}\to \M_{g,k}$. The latter is defined by forgetting the map to $X$ and the last $k+n$ marked points, and contracting those components of the curve which have become unstable. We follow the exposition in Appendix 2 of \cite{CGi} to relate descendent and ancestor potentials. The geometry of this relationship goes back to the paper of Kontsevich-Manin \cite{KM2} and Getzler \cite{Ge}.

Let $L$ be one of the cotangent line bundles over $X_{g,k+l+n,d}$ (say, the 1st one), and $\bar{L}$ its counterpart pulled back from $\M_{g,k}$. They are the same outside the locus where the $1$st marked point lies on a component to be contracted. This shows that there is a holomorphic section of $Hom (\bar{L}, L)$ vanishing on the virtual divisor $j: D\to X_{g,n,d}$ formed by gluing genus $g$ stable maps, carrying all but the 1st out of the first $k$ marked points, with genus $0$ stable maps, carrying the 1st one. In fact, like in the case of the of WDVV-equation, the divisor has self-intersections (see Figure 2), and we have to refer once again to \cite{GiK} for a detailed discussion of the K-theoretic exclusion-inclusion formula
\[ \O - \O(-D)= j_*\O_D-j_*\O_{D_{(2)}}+j_*\O_{D_{(3)}}-\cdots \]
which expresses $1-\bar{L}/L = \O-\O(-D)$ in terms of structure sheaves of the strata $D_{(m)}$ of $m$-tuple self-intersections.   

\begin{figure}[htb]
\begin{center}
\epsfig{file=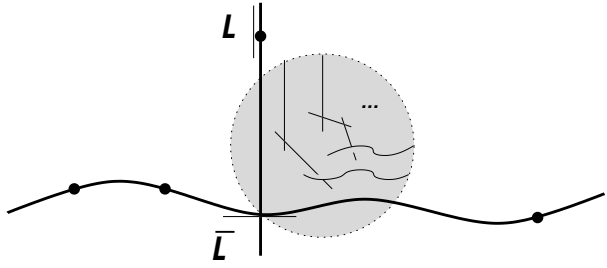} 
\end{center}
\caption{The divisor $D$ and its self-intersections}
\end{figure}

We will use this relationship to rid systematically of $L$'s in favor of $\bar{L}$'s in the correlators. For this, we will have to consider
``mixed'' correlators, which allow both $L$ and $\bar{L}$ at the same seat.
Let us use the notation $\lan \phi L^a\bar{L}^b|$ in correlator expressions which have the specified inputs (here $\phi \in K^0(X)\otimes \gL$) in the singled out (first) seat, provided that all other inputs in all terms of the expression are the same. For $a>0$, we have:
\begin{align*} \lan \phi L^a\bar{L}^b| =& \lan \phi L^{a-1}\bar{L}^{b+1}| +  \lan \phi L^a\bar{L}^b (1-\bar{L}/L) | \\
= &\lan \phi L^{a-1}\bar{L}^{b+1}| + \sum_{\a,\b}\llan \phi L^a, \phi_{\a}\rran_{0,2}G^{\a\b} \lan \phi_\b \bar{L}^b|\ .\end{align*}
Note that marked points (not shown on Figure 2) which carry the inputs $\tau$ or permutable inputs $t$ can be distributed in any way between the components of the curve, and the above factorization of correlators under gluing is justified by the {\em permutation-equivariant binomial formula} from Part I.  

Iterating the procedure, we have:
\begin{align*} \lan \phi L^a| = \sum_{\a,\b} &\llan \phi L^a, \phi_{\a}\rran_{0,2}G^{\a\b} \lan \phi_{\b}|  + \cdots + \sum_{\a,\b}\llan \phi L, \phi_{\a}\rran_{0,2}G^{\a\b}\lan \phi_{\b}\bar{L}^{a-1}|\\ &+\lan \phi \bar{L}^a| 
  = \lan \phi \bar{L}^a|+\sum_{b=0}^{a-1}\sum_{\a,\b}\llan \phi L^{a-b},\phi_{\a} \rran_{0,2}G^{\a\b} \lan \phi_{\b} \bar{L}^b| .\end{align*}  
Similarly, for negative exponents, we have  
\begin{align*} \lan \phi L^{-a-1}|&= \lan \phi L^{-a}\bar{L}^{-1}|-\lan (1-\bar{L}/L)L^{-1}/\bar{L}| \\ &= \lan \phi L^{-a}\bar{L}^{-1}| - \sum_{\a,\b}\llan \phi L^{-a},\phi_{\a}\rran_{0,2}G^{\a\b} \lan \phi_{\b} \bar{L}^{-1}| \\ & =\lan \phi L^{-a+1}\bar{L}^{-2}| - \cdots \\
 = \lan \phi \bar{L}^{-a-1}| &-\sum_{b=0}^{a}\sum_{\a,\b}\llan \phi L^{-a+b}, \phi_{\a}\rran_{0,2}G^{\a\b} \lan \phi_{\b} \bar{L}^{-b-1}| . \end{align*}
In fact the result can be concisely described as
\[ \lan \x(L) | = \left\lan \left[ S(\bar{L})\, \x(\bar{L})\right]_{+} \middle|, \right.\]
where the operator $S$ is as in the previous section: 
\[ S (q)\, \phi = \sum_{\a,\b}\left((\phi,\phi_{\a})+\llan \frac{\phi}{1-L/q}, \phi_{\a} \rran_{0,2} \right) G^{\a\b}\phi_{\b},\]
and $[f(q)]_{+}$ means extracting from a rational function of $q$ the Laurent polynomial part. The latter procedure, understood as projection along the space of rational functions regular at $q=0$ and vanishing at $q=\infty$, can be described by Cauchy's residue formula:
\[ [f(q)]_{+} = -\Res_{w=0,\infty} \frac{f(w)dw}{w-q}.\]
For a Laurent polynomial $\x(q)$ we have:
\[ -\Res_{w=0,\infty} \frac{\x(w)}{(1-L/w)}\, \frac{dw}{(w-\bar{L})} =
\frac{\x(L)}{1-\bar{L}/L}-\frac{\x(\bar{L})}{1-L/\bar{L}} .\]
For $\x(L)=L^a$ we get $(L^{a+1}-\bar{L}^{a+1})/(L-\bar{L})=\sum_{b=0}^{a}L^{a-b}\bar{L}^b$,
and hence
\begin{align*} \left\lan \left[ S(L)\, \phi \bar{L^a}\right]_{+} \middle| \right. & = \\
\sum_{\a,\b} & \left( (\phi,\phi_{\a})G^{\a\b}\lan \phi_{\b} \bar{L}^a| + \sum_{b=0}^{a}\llan \phi L^b,\phi_{\b} \rran_{0,2}G^{\a\b}\lan \phi_{\b} \bar{L}^{a-b}| \right) ,\end{align*} 
which agrees with what we found earlier, because
\[ \sum_{\a}(\phi,\phi_{\a})G^{\a\b}=(\phi,\phi_{\b})-\sum_{\a}\llan \phi,\phi_{\a}\rran_{0,2}G^{\a\b}.\]
For $x(L)=L^{-a-1}$, it works out similarly:
\[ \frac{L^{-a-1}}{1-\bar{L}/L}+\frac{\bar{L}^{-a-1})}{1-L/\bar{L}} =
=\frac{L^{-a}-\bar{L}^{-a}}{\bar{L}^{-1}-L^{-1}}
=-\sum_{b=1}^{a}L^{-a+b}\bar{L}^{-b+1}.\]

The same procedure can be applied to each seat in the correlators. We conclude that for stable values of $(g, m)$,
\[ \llan \x(L),\dots,\x(L) \rran_{g,m} = \llan \y(\bar{L}),\dots,\y(\bar{L})\rran_{g,m},\ \text{where}\ \y(q)=\left[ S(q)\x(q)\right]_{+}\]
quite analogously to the cohomological results of \cite{KM2}.  Here all the correlators, as well as $S$, depend on the permutable parameters $t$ and non-permutable $\tau$. For a fixed $t$, assembling the correlators into the generating functions $\F_g$ and $\bar{\F}_g$, we find:
\[ \F_g(\tau+\x)=\bar{\F}_g^{(\tau)}([S_{\tau}\x]_{+}) + \delta_{g,1} \llan\ \rran_{1,0}^{(\tau)}+\delta_{g,0}\sum_{m=0}^2\frac{1}{m!}\llan \dots, \x(L), \dots \rran_{0,m}^{(\tau)},\]
where the decorations by $\tau$ remind on the dependence on the parameter, and the terms on the right represent correlators with unstable values of $(g,m)=(1,0),(0,0),(0,1),(0,2)$, present in descendent, but absent in ancestor potentials.

Now we engage the shift of the origin $\x = \y+q-1-t-\tau$. We have:
\begin{align*} \left[\frac{q-1-t-\tau}{1-L/q}\right]_{+} & = -\Res_{q=0,\infty} \frac{w-1-t-\tau}{1-L/w}\frac{dw}{w-q} =\\
  \frac{q-1-t-\tau}{1-L/q} &+\frac{L-1-t-\tau}{1-q/L} = L+q-1-t-\tau.\end{align*}
Therefore $[S(q) (q-1-t-\tau) ]_{+} = $
\begin{align*} &\sum_{\a,\b}\left((q-1-t-\tau, \phi_{\a})+\llan \left[\frac{q-1-t-\tau}{1-L/q}\right]_{+}, \phi_{\a}\rran_{0,2}\right) G^{|a\b}\phi_{\b}\\ &= q-1-t-\tau+\sum_{\a,\b}\llan L, \phi_{\a}\rran_{0,2}G^{\a\b}\phi_{\b} = q-1. \end{align*}
The last equality is due to the string and dilaton equations:
\begin{align*}  \llan L, \phi_{\a}\rran_{0,2}&=\sum_{d,l,n}\frac{Q^d}{l!}\lan L, \phi_{\a}, \tau,\dots, \tau; t,\dots, t\ran^{S_n}_{0,2+l+n,d}= \\
  \lan L, \phi_{\a}, \tau+t\ran_{0,3,0}&+\sum_{d,l,n}\frac{Q^d}{l!}\lan \phi_{\a}, \tau+t, \tau,\dots,\tau; t,\dots,t\ran_{0,2+l+n}^{S_n} = \\
  (\tau+t ,\phi_{\a}) &+ \llan \tau+t, \phi_{\a}\rran_{0,2}, \end{align*} 
and hence $\sum_{\a,\b}\llan L, \phi_{\a}\rran_{0,2}G^{\a\b}\phi_{\b}=t+\tau$.

Finally, using the dilaton equations
\begin{align*} \llan L-1, A \rran_{0,2}&=-\llan A\rran_{0,1}+\llan A, t+\tau\rran_{0,2}, \\ \llan L-1\rran_{0,1}&=-2\llan \ \rran_{0,0}+\llan t+\tau\rran_{0,1},\end{align*}
we find that 
\[ \llan \ \rran_{0,0}+\llan \y+L-1-t-\tau \rran_{0,1}+\frac{1}{2}\llan\y+L-1-t-\tau ,\y+L-1-t-\tau\rran_{0,2}\]
transforms into $\llan \y,\y\rran_{0,2}/2$. Indeed, the terms linear in $\y$
\[ \llan \y\rran_{0,1}+\llan L-1-t-\tau,\y\rran_{0,2}=\llan \y\rran_{0,1}-\llan \y\rran_{0,1}+\llan \y, t+\tau\rran_{0,2}-\llan t+\tau,\y\rran_{0,2}\]
cancel out. The $\y$-independent terms
\begin{align*} &\frac{1}{2}\llan L-1,L-1\rran_{0,2}-\llan L-1, t+\tau\rran_{0,2}+\frac{1}{2}\llan t+\tau,t+\tau\rran_{0,2}\\
  &+\llan L-1\rran_{0,1}-\llan t+\tau\rran_{0,1}+\llan \ \rran_{0,0}= -\frac{1}{2}\llan L-1,t+\tau\rran_{0,2}+ \\ &\frac{1}{2}\llan t+\tau,t+\tau\rran_{0,2}+\frac{1}{2}\llan L-1\rran_{0,1}-\llan t+\tau\rran_{0,1}+\llan \ \rran_{0,0}= \\ &\frac{1}{2}\llan t+\tau\rran_{0,1}
+\frac{1}{2}\llan t+\tau\rran_{0,1}-\llan \ \rran_{0,0}-\llan t+\tau\rran_{0,1}+\llan \ \rran_{0,0} . 
\end{align*}
cancel out too. Thus, we obtain
\[ \F_g(\y+q-1-t)=\bar{\F}_g([S \y]_{+}+q-1)+\delta_{g,1}\llan \ \rran_{1,0}+ \frac{\delta_{g,0}}{2}\llan \y(L),\y(L)\rran_{0,2}.\]
In view of Proposition 5 from the previous section, we have proved the following theorem.

\medskip

{\tt Theorem 1.} {\em The total descendent potential after the shift by $1-q+t$:
  \[ \D (1-q+t+\x) = e^{\sum_{g\geq 0} \h^{g-1}\F_g(\x)},\]
and the $\tau$-family of total ancestor potentials after the shift by $1-q$:
  \[ \A_{\tau} (1-q+\x):=e^{\sum_{g\geq 0} \h^{g-1}\bar{\F}_g^{(\tau)}(\x)}, \ \tau\in K^0(X)\otimes \gL,\]
are related by the family of quantized operators
  \[ \D = e^{F_1(\tau)}\,\widehat{S^{-1}_{\tau}}\, \A_{\tau},\]
  where $F_1(\tau)=\llan\ \rran_{1,0}:=\sum_{k,n,d}\frac{Q^d}{k!}\lan \tau,\dots, \tau;t,\dots,t\ran_{1,k+n,d}^{S_n}$ is the generating function for primary GW-invariants of genus 1.}

\medskip

Passing to the quasi-classical limit $\h\to 0$, one obtains

\medskip

{\tt Corollary 1.} {\em The graph $\L \subset \K$ of the differential of the genus-0 descendent potential $\F_0(\y+q-1-t)$ and the $\tau$-family $\L^{(\tau)} \subset \K$ of the graphs of the differentials of genus-0 ancestor potentials $\bar{F}_0^{(\tau)}(\y+q-1)$ are related by symplectic transformations $S_{\tau}: (\K, \Omega) \to (\K, \bar{\Omega}^{(\tau)})$: 
\[ \L = S_{\tau}^{-1}\L^{(\tau)}.\]} 

The genus-0 ancestor correlators $\llan \x(\bar{L}),\dots,\x(\bar{L})\rran_{0,m}$ have the ``zero 2-get'' property \cite{Ge2}; namely they have zero 2-jet along the subspace $\x \in \K_{+}$, where $\x(1)=0$. This is because $\bar{L}_i$ are pull-backs of the line bundles $L_i$ from the Deligne-Mumford space $\M_{0,m}$, which is a {\em manifold} of dimension $m-3$, and where therefore any product of $m-2$ factors $\bar{L}_i-1$ vanishes for dimensional reasons. Since the dilaton shift $1-q$ also vanishes at $q=1$, we conclude that {\em $\L^{|tau}$ is tangent to $\K_{+}$ along $(1-q)\K_{+}$. Consequently, $T_{|tau}:= S_{\tau}^{-1} \K_{+}$ is tangent to $\L$
along $(1-q)T_{\tau} \subset \L$.} In fact, as $\tau$ varies, these spaces sweep $\L$ (which is easy to check modulo Novikov's variables, and then apply the formal Implicit Function Theorem.)  
   
\medskip

{\tt Corollary 2.} {\em $\L \subset (\K, \Omega)$ is an overruled Lagrangian cone, i.e. its tangent spaces $T:=T_{\J}\L$ are tangent to $\L$ exactly along $(1-q)T\subset \L$.}

\section*{Example: $X=pt$}

In Part I, we found that for $t\in \gL$
\[\J(0,t):=1-q+t+\sum_{n\geq 2} \lan\frac{1}{1-qL};t,\dots,t\ran_{0,1+n}^{S_n}=(1-q)e^{\sum_{k>0}\Psi^k(t)/k(1-q^k)}.\]
It follows from the string equation (see Corollary 2 of Proposition 4) that for $\tau \in \gL$
\[ \J(\tau,t)=1-q+t+\tau+\llan \frac{1}{1-qL}\rran_{0,1}=(1-q)e^{\tau/(1-q)+\sum_{k>0}\Psi^k(t)/k(1-q^k)}.\]
Taking $q=0$ (and using the string equation twice), we find the variable metric
\[ G(\tau)=G_{11}(\tau) := \llan 1, 1, 1\rran_{0,3}=1+\tau+t+\llan 1 \rran_{0,1}=e^{\tau+\sum_{k>0} \Psi^k(t)/k}.\]
Using the string equation once more, we derive that
\begin{align*} S_{\tau}(q) &:= \left(1+\llan \frac{1}{1-L/q},1\rran_{0,2}\right)G^{-1}(\tau)=
  \frac{\J (1/q)}{1-1/q}\, G^{-1}(\tau)\\ &= e^{\textstyle \tau/(q-1)-\sum_{k>0}\Psi^k(t)/k(q^k-1)},
  \\  S^{-1}_{\tau}(q)&=e^{\textstyle \tau/(1-q)+\sum_{k>0}\Psi^k(t)/k(1-q^k)} =\frac{\J(q)}{1-q}, \end{align*}
and find the range $\L_t$ of the J-function $\K_{+}\to \K:\x\mapsto \J(\x, t)$ to be
\[ \L_t = \bigcup_{\tau\in \gL}\, e^{\textstyle \tau/(1-q)+\sum_{k>0}\Psi^k(t)/k(1-q^k)} (1-q)\K_{+}. \]
At $\tau=0$, we have here one of the subspace in $\K$, depending on $t$, whose union over $t\in \gL_{+}$, according to the results of Part III, yields the range $\L$ of the permutation-equivariant J-function $\t \mapsto \J(0,\t)$
\[  \L=\bigcup_{t\in \gL_{+}}\, e^{\textstyle \sum_{k>0} \Psi^k(t)/k(1-q^k)} (1-q)\K_{+}.\]
In fact this picture remains true in general, as we will now show.
  
\section*{Adelic characterization}

We return now to the mixed genus-0 descendent potential $\F_0(\x,\t)$ with the permutable input $\t \in \K_{+}$ allowed to involve the cotangent line bundles $L_i$. In the symplectic loop space $(\K,\Omega)$, it is represented by the dilaton-shifted graph of its differential. According to Proposition 4 and its Corollary 1, it is the range of the J-function 
\begin{align*} \K_{+}\ni \x &\mapsto \J(\x,\t):=1-q+\t(q)+\x(q)+\\
  &\sum_{\a,k,n,d}\phi^{\a}\frac{Q^d}{k!}\lan \frac{\phi_{\a}}{1-qL}, \x(L),\dots,\x(L);\t(L),\dots,\t(L)\ran_{0,1+k+n,d}^{S_n},\end{align*}
and has the form of a Lagrangian cone $\L_{\t}$, depending on the parameter $\t\in \K_{+}$. According to the results of the previous section $\L_{\t}$ is an overruled Lagrangian cone whenever $\t$ is constant in $q$.   
We combine this information with the {\em adelic characterization} of the J-function given in \cite{GiTo, ToK, ToT}\footnote{Formally speaking, there only the case $\t=0$ is considered, but the results extend without change to the general case, where the moduli orbi-spaces are $X_{0,1+k+n,d}/S_n$ rather than $X_{0,1+k,d}$.}  and discussed in Part III, to prove the following theorem.

\medskip

{\tt Theorem 2.} {\em The range $\L$ of permutation-equivariant J-function
  $\t \mapsto \J(0,\t)$ (with $\t \in \K_{+}$, and $\t(1)\in K^0(X)\otimes \gL_{+}$, where $\gL_{+}$ is a certain neighborhood of $0\in \gL$) has the form 
  \[ \L = \bigcup_{t\in K^0(X)\otimes \gL_{+}} \, (1-q)\, S_0^{-1}(q)_t\,\K_{+},\]
  where the operators $S_{\tau}(q)$ evaluated at $\tau=0$ still depend on the parameter $t\in K^0(X)\otimes \gL_{+}$:
  \[ S_0^{-1}(q)_t\, \psi := \psi + \sum_{\a} \phi_{\a} \sum_{n,d} Q^d \lan \psi, \frac{\phi_{\a}}{1-qL} ; t,\dots, t\ran_{0,2+n,d}^{S_n}.\]} 

\medskip

{\tt Proof.} According to the adelic characterization results, {\em a rational function $\f \in \K$ lies in $\L_{\t}$ if and only if its Laurent series expansions $\f_{(\z)}$ near $q=1/\z$ satisfy the following three conditions:

  (i) $\f_{(1)} \in \L^{fake}\subset \hat{\K}$, the range, in the space $\hat{\K}$ of Laurent series in $q-1$ with vector coefficients in  $\K^0(X)\otimes \gL$, of the J-function in the {\em fake} quantum K-theory of $X$;

  (ii) when $\z\neq 0,1,\infty$ is a primitive $m$th root of unity,
  $\f_{(\z)}(q^{1/m}/\z) \in \L^{(\z)}_{\t}$, a certain Lagrangian subspace in $\hat{K}$ which will be specified below;

  (iii) when $\z \neq 0,\infty$ is not a root of unity, $\f_{(\z)}$ is a power series in $q-1/\z$, i.e. $\f$ has no pole at $q=1/\z$.}

\medskip

To elucidate the situation, recall that in {\em fake} K-theory, the genuine holomorphic Euler characteristics $\chi ({\mathcal M} ; V)$ are replaced with their ``fake'' values given by the right-hand-side of the Hirzebruch--Riemann--Roch formula:
\[ \chi^{fake}({\mathcal M}; V) := \int_{[{\mathcal M}]} \ch (V) \td ({T_{\mathcal M}}).\]
Fake in this sense GW-invariants were studied, e.g. in \cite{Co}. In particular, the range of the fake J-function is known to be an overruled Lagrangian cone $\L^{\fake}\subset (\hat{\K}, \hat{\Omega})$, where
$\hat{\Omega}(\f,\g) = \Res_{q=1}(\f(q^{-1}),\g(q))\, q^{-1}dq$.

The moduli spaces of stable maps behave as virtual {\em orbifolds} (rather than manifolds), and the genuine holomorphic Euler characteristics are given by the virtual Kawasaki--Riemann--Roch formula \cite{ToK}, summing up certain fake holomorphic Euler characteristics of the inertia orbifold (of the moduli spaces $X_{0,1+k+n}/S_n$ in our situation). Figure 3, essentially copied from Part III, is to remind us of the recursive device keeping track of all Kawasaki contributions into the J-function.

 %\left( S_{\tau (\J_X(\t)_{(1)}}) \right) \hat{\K}^X_{+}

\begin{figure}[htb]
\begin{center}
\epsfig{file=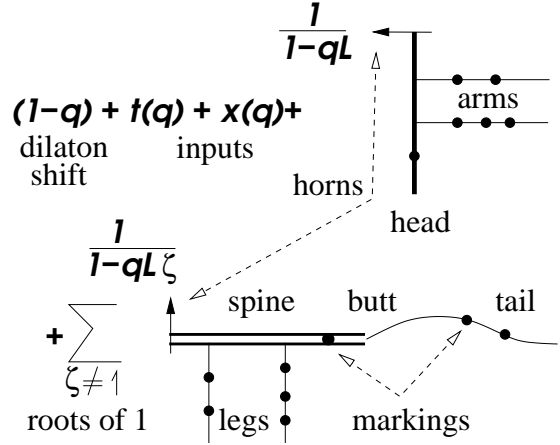} 
\end{center}
\caption{Adelic characterization}
\end{figure}

In particular, it shows that the values of the J-function, when expanded into near $q=1$, lie in $\L^{fake}$, and when expanded near a primitive $m$th root of $1$, they are characterized in terms of certain twisted fake invariants of the orbifold target space $X\times B\ZZ_m$. The latter, in their turn, are expressed in terms of the untwisted fake invariants of $X$. Namely in the test (ii) above, the subspace $\L_{\t}^{(\z)} \subset \hat{K}$ is obtained from {\em a certain tangent space $T^{fake}_{\t}$ to $\L^{fake}$ by the linear transformation:
  \[ \L_{\t}^{(\z)} = e^{\textstyle \sum_{k>0}\left(  \frac{\Psi^k(T^*_X)}{k(1-\zeta^{-k}q^{k/m})}-\frac{\Psi^{km}(T_X^*)}{k(1-q^{km})}\right)}
  \Psi^m (T^{fake}_{\t})\otimes_{\Psi^m(\gL)} {\gL}. \]
}
In our present discussion, it is important to figure out what determines the application point of the tangent space $T_{\t}^{fake}$. On the diagram, it is determined by {\em legs}, which are related by the Adams operation $\Psi^m$ to {\em arms} (see Part III, or \cite{GiTo}). Note however, that the markings on the legs (each representing $m$ copies of markings on the arms attached to the $m$-fold cover of the spine curve) are allowed to carry permutable inputs $\t$, but not allowed to carry the non-permutable inputs $\x$, because their numbering would break the $\ZZ_m$-symmetry of the covering curve). Consequently, {\em the value of $\J^{fake} \in \L^{fake}$, which determines the application point of the tangent space $T_{\t}$, is obtained by the expansion near $q=1$ of the J-function with the non-permutable input $\x= 0$:
  $T_{\t} = T_{\J(0,\t)_{(1)}}\L^{fake}$. } In fact, since $\L^{fake}$ is overruled, its tangent spaces to $\L^{fake}$ are parameterized by $K^0(X)\otimes \gL$. Let us analyze the map $\t \mapsto\ \text{(tangent space to $\L^{fake})$}$.

In degree $d=0$, the J-function of $X$ coincides with the J-function of the point target space with coefficients in the $\lambda$-algebra $\gL' := K^0(X)\otimes \gL$. It was described in section Example. For $\t = t\in \gL'$
(i.e. $q$-independent), we have
\[ \J(0, t)_{(1)} = (1-q)e^{\textstyle \sum_{k>0} \Psi^k(t)/k^2(1-q)}\times (\text{power series in $q-1$}).\]
In other words, $\sum_{k>0} \Psi^k(t)/k^2$ is the parameter value of the tangent
space to $\L^{fake}$ associated to the input $t\in \gL'$ in this approximation.

The series is not guaranteed to converge. E.g. under the identification of $K^0(X)\otimes \QQ$ with $H^{even}(X,\QQ)$ by the Chern character, $\Psi^k$ acts on $H^{2r}(X)$ as multiplication by $k^r$, and the series $\sum_{k>0} k^{r-2}$ diverges unless $r=0$. To handle this difficulty, we assume that the ground ring $\gL$ is topologized with a filtration $\gL\supset \gL_{+}\supset \gL_{++}\supset \dots $ by ideals such that $\Psi^k$ with $k>1$ increase the filtration. For instance, when $\gL=\QQ[[Q]]$ is the Novikov ring, $\Psi^k(Q^d)=Q^{kd}$, the filtration by the powers of the maximal ideal is taken. When $\gL=\QQ [[N_1,N_2,\dots]]$ is the ring of symmetric functions, $\Psi^k(N_r)=N_{kr}$, the filtration by degrees of symmetric functions suffices. 
Then the map $t\mapsto \sum_{k>0} \Psi^k(t)/k^2$ converges for $t\in K^0(X)\otimes \gL_{+}$, and is invertible in this range,\footnote{Even in the entire $H^0(X,\gL)$, if $\sum_{k>0}k^{-2}=\pi^2/6$ is adjoined to $\gL$.} since $\Psi^1(t)=t$.

Returning to the general input $\t$ and degree $d\geq$, we conclude from the formal Implicit Function Theorem, that there is a well-defined map
\[ \T: \{ \t \in \K_{+}\, |\, \t(1) \in K^0(X)\otimes \gL_{+} \} \to  K^0(X)\otimes \gL_{+}, \]
such that $T_{\J(0,\t)_{(1)}} \L^{fake}= T_{\J(0,\T(\t))_{(1)}}\L^{fake}$. For all inputs $\t$ with the same value $\T(\t)$, the adelic characterization tests (i), (ii), (iii) coincide.    

By the same token, for each $t$ there is a well-defined map $\K_{+} \to K^0(X)\otimes \gL: \x \mapsto \tau(\x)$, such that $T_{\J(\x,t)}\L_t = T_{\J(\tau(\x),t)}\L_t$. For all inputs $\x$ with the same $\tau(\x)$, the values $\J (\x,t)$ of the J-function form the ruling space
$(1-q)\, S_{\tau}(q)_t\, \K_{+}$ of the overruled cone $\L_t$. For all such points, the localizations $\J(\x,\t)_{(1)}$ lie in the same ruling space of $\L^{fake}$, and moreover, when $\tau =0$, the last ruling space is the one where $\J(0,\t)$ with $\T(\t)=t$ lie.  Thus, for rational functions from the space $(1-q)\, S_0(q)_t\, \K_{+}$ and for the values $\J(0,\t)$ with $\T(\t)=t$, the adelic characterization tests (i), (ii), (iii) coincide, i.e.
the localizations in test (i) lie in the same ruling space of $\L^{fake}$, and the tangent spaces to $\L^{fake}$ involved into test (ii) are the same. 
Therefore the two sets of rational functions coincide:
\[ \{ \J(0,\t) \, | \, \T(\t)=t\} = (1-q)\, S_0(q)_t\, \K_{+} .\]  
Taking the union over $t\in \gL_{+}$ completes the proof. \qed

\medskip

{\tt Corollary 1.} {\em $\L_{\t} = \L_t$, where $t=\T(\t)$.}

\medskip

{\tt Corollary 2.} {\em Each $\L_{\t}$ is an overruled Lagrangian cone invariant under the string flow $\f \mapsto e^{\epsilon /(1-q)}\f$, $\epsilon\in \gL$.}  

\medskip  

{\tt Remark.} The range $\L \subset (\K, \Omega)$ of the permutation-equivariant J-function $\t \mapsto \J(0,\t)$ is a cone ruled by the family $t\mapsto R_t:=(1-q)S_0(q)_t\K_{+}$ of {\em isotropic} subspaces (and is in this sense ``overruled'') but it is not Lagrangian, nor is it invariant under the string flow, as the example of $X=pt$ readily illustrates.
In particular, the spaces  $R_t/(1-q)$ are not tangent to $\L$, and do not form
semi-infinite variations of Hodge structures in the sense of S. Barannikov \cite{Bar}. Nevertheless from Proposition 2 (dilaton equation), we have:

\medskip

{\tt Corollary 3.} {\em The permutation-equivariant genus-0 descendent potential
\[ \F_0(0,\t):=\sum_{0,n,d}Q^d\lan \t(L),\dots,\t(L)\ran_{0,n,d}^{S_n} \]
is reconstructed from the permutation-equivariant J-function by
\[ \frac{1}{2}\Omega \left( [\J (0,\t)]_{-}, [\J (0,\t)]_{+}\right) = \F_0(0,\t) + \frac{(\Psi^2(\t(1)), 1)}{2}.\]}

\enddocument